\documentclass[9pt]{article}
\usepackage{amsmath}
\usepackage{amsfonts}
\usepackage{amssymb}

\newtheorem{lemma}{Lemma}[section]

\newtheorem{theorem}{Theorem}[section]

\let\Section=\section
\def\section{\setcounter{equation}{0}\Section}
\sf

\begin{document}
\title{Multiple solutions for a nonhomogeneous Schr\"odinger-Maxwell system in $\mathbb{R}^3$ }
 \author{ \small Yongsheng  Jiang,  Zhengping  Wang and  Huan-Song Zhou \thanks{
 \indent Corresponding author: hszhou@wipm.ac.cn; Phone: +86-27-87199196; Fax: +86-27-87199291. This work was supported by NSFC(10801132, 11071245 and 11126313)}\\
 \\ \small Wuhan Institute of Physics and Mathematics,  Chinese Academy of
Sciences \\ \small P.O. Box 71010, Wuhan 430071, P.R. China}
 \date{}
\maketitle

\noindent {\bf Abstract:}  The paper considers the  following  nonhomogeneous Schr\"odinger-Maxwell
system
\begin{equation}\nonumber
\left\{\begin{array}{ll}
 -\Delta u +  u+\lambda\phi (x) u =|u|^{p-1}u+g(x),\ x\in \mathbb{R}^3, \\
 -\Delta\phi = u^2, \ x\in \mathbb{R}^3,
 \end{array}\right. \leqno{(SM)}
\end{equation}
where $\lambda>0$, $p\in(1,5)$ and $g(x)=g(|x|)\in
L^2(\mathbb{R}^3)\setminus\{0\}$.
 There seems no any results on the existence of multiple solutions to problem (SM) for $p \in (1,3]$.
 In this paper, we find that there is a constant$C_p>0$ such that problem (SM)  has at least two solutions for all $p\in (1,5)$ provided
 $\|g\|_{L^2} \leq C_p$, but only for $p\in(1,2]$ we need $\lambda>0$ is small. Moreover, $C_p=\frac{(p-1)}{2p}[\frac{(p+1)S^{p+1}}{2p}]^{1/(p-1)}$, where $S$ is the Sobolev constant.

\noindent {\it Keywords}: Multiple
solutions, nonhomogeneous,  Schr\"odinger-Maxwell system. \\
\noindent {\it 2010 Mathematical Subject Classification}: 35J20,
35J60.

\section{Introduction}
In this paper, we are concerned with the existence of multiple
solutions for the following nonhomogeneous Schr\"odinger-Maxwell
system
\begin{equation}\label{eq:1.1}
\left\{\begin{array}{ll}
 -\Delta u +  u+\lambda\phi (x) u =|u|^{p-1}u+g(x),\ x\in \mathbb{R}^3, \\
 -\Delta\phi = u^2, \ x\in \mathbb{R}^3,
 \end{array}\right.
\end{equation}
where $\lambda>0$ is a parameter, $p\in(1,5)$ and $g(x)=g(|x|)\in
L^2(\mathbb{R}^3)$.

Problem (\ref{eq:1.1}) is related to the study of nonlinear Schr\"odinger
equation for a particle in an electromagnetic field. For more details on the physical aspects about the problem we refer the
reader to \cite{VBenciDFortucto-TopMNA,SalvatoreA-ANonSt} and the
references therein. If
$g(x)\equiv0$, the existence of solutions to problem (\ref{eq:1.1})
has been discussed under different ranges of $p$, for examples,
\cite{GMCoclite-CAA} \cite{D'AprileTMugnaiD-PRSESM}
\cite{d'AveniaP-ANonSt} for $p\in[3,5)$,
\cite{AzzolliniAandPomponioA-Jmaa} for $p\in(2,5)$,
\cite{KikuchiHiroaki-NonlAnal} for $p\in[2,3)$ and
\cite{AmbrosettiAntonioRuizDavid-CCM}
\cite{AzzolliniAveniaPomponio-AIHP} \cite{RuizD-JFA} for $p\in(1,5)$
or general nonlinearity, etc. \iffalse For some $\lambda_0>0$, it
was proved in \cite{KikuchiHiroaki-NonlAnal} that (\ref{eq:1.1})
with $\mu=0$ and $p\in(1,2]$ has no nontrivial solution if
$\lambda\geqslant\lambda_0$. Ruiz in \cite{RuizD-JFA} proved that
$\lambda_0=\frac{1}{4}$ and that (\ref{eq:1.1}) has solution if
$\lambda>0$ small enough.\fi Some recent results in this direction
was summarized in \cite{AmbrosettiA-MJM}. However, if
$g(x)\not\equiv 0$, only a few results are known for problem
(\ref{eq:1.1}) when $p\in(3,5)$.
% and with some special potential functions,
%such as \cite{SJChen-CLTang,SalvatoreA-ANonSt,MBYang-BRLi}.
In \cite{SalvatoreA-ANonSt}, three radially symmetric solutions of
(\ref{eq:1.1}) were obtained for  $p\in(3,5)$ and $\|g\|_{L^2}$ is
small enough. In \cite{SJChen-CLTang,MBYang-BRLi}, the authors
considered problem (\ref{eq:1.1}) with certain potential and the
existence of multiple solutions is established for $p\in(3,5)$. To
the authors' knowledge, it is still open whether the problem
(\ref{eq:1.1}) has multiple solutions under $p\in(1,3]$ and
$g(x)\not\equiv0$. The aim of this paper is to prove that problem
(\ref{eq:1.1}) has at least two solutions for all $p\in(1,5)$ and
$\|g\|_{L^2}$ is suitably small.

For $u\in H^1(\mathbb{R}^3)$, let $\phi_u$ be the unique solution of
$-\Delta \phi=u^2$ in $D^{1,2}(\mathbb{R}^3)$, then
\begin{equation}\label{eq:1.2}
\phi_u(x)=\frac{1}{4\pi}\int\frac{u^2(y)}{|x-y|}dy,
\end{equation}
here and in what follows, we denote $\int_{\mathbb{R}^3}$ simply by $\int$. Define the energy functional $I_\lambda:
H^1(\mathbb{R}^3)\to \mathbb{R}$ by
\begin{equation}\label{eq:1.3}
I_\lambda(u)=\frac{1}{2}\int|\nabla u|^2+ u^2
dx+\frac{\lambda}{4}\int\phi_uu^2dx-\frac{1}{p+1}\int|u|^{p+1}dx-\int
g(x)udx.
\end{equation}
 If $g(x)\in
L^2(\mathbb{R}^3)$ and $p\in [1,5]$, it is known that $I\in
C^1(H^1(\mathbb{R}^3),\mathbb{R})$ and for any $\varphi\in
H^1(\mathbb{R}^3)$,
\begin{equation}\label{eq:1.4}
\langle I_\lambda'(u),\varphi\rangle=\int\nabla u\nabla \varphi+
u\varphi dx+\lambda\int\phi_u u\varphi dx-\int|u|^{p-1}u\varphi
dx-\int g(x)\varphi dx.
\end{equation}
Furthermore, if $g(x)=g(|x|)$ and $u\in H_r^1(\mathbb{R}^3)$
satisfies $I_\lambda'(u)\varphi=0$ for  all $\varphi\in
H_r^1(\mathbb{R}^3)$, Lemma 2.4 of \cite{D'AprileTMugnaiD-PRSESM}
showed that $(u,\phi_u)$
 satisfies (\ref{eq:1.1}) in the weak sense. For simplicity, in many cases
we just say $u\in H_r^1(\mathbb{R}^3)$, instead of $(u,\phi_u)\in
H_r^1(\mathbb{R}^3)\times D^{1,2}(\mathbb{R}^3)$, is a weak
solution of (\ref{eq:1.1}).\\

For all $\lambda>0$, $p\in (1,5)$ and  $\|g\|_{L^2}$ suitably small,
it is not difficult to get a solution $u_0$ of (\ref{eq:1.1}) by the
Ekeland's variational principle. Moreover, $u_0$ is a local
minimizer of $I_\lambda$ and of negative energy, that is, $I_\lambda(u_0)<0$. To get a solution
of (\ref{eq:1.1}) with positive energy, we have to study the problem
(\ref{eq:1.1}) in the following two cases: $p\in (2,5)$ and $p\in
(1,2]$, respectively.

When $p\in (2,5)$, by using the transform $w_t(x)=t^2w(tx)$ for some
$w\in H_r^1(\mathbb{R}^3)$ and $t>0$ large enough, we can show that
$I_\lambda$ satisfies the mountain pass geometry for any $\lambda>0$
(see Lemma \ref{Le3.1}) and get a $(PS)_c(c>0)$ sequence $\{u_n\}$
of $I_\lambda$. For $p\in [3,5)$, it is easy to prove the
boundedness of $\{u_n\}$ and the (PS) condition. But for $p\in
(2,3)$, it is still not clear if the (PS) condition holds. To
overcome this difficulty in the case of $g(x)\equiv 0$, Ruiz
\cite{RuizD-JFA} introduced an interesting manifold $\mathcal{M}$
and then proved that there exists a positive radial function
$\tilde{u}$ such that
$0<I_\lambda(\tilde{u})=\inf\{I_\lambda(u):u\in\mathcal{M}\}$ and
$I_\lambda'(\tilde{u})=0$. Using this manifold $\mathcal{M}$ and the
concentration compactness principle,  Azzollini and Pomponio
\cite{AzzolliniAandPomponioA-Jmaa} established the existence of a
ground state for problem (\ref{eq:1.1}) under $g(x)\equiv 0$ and
$p\in (2,5)$. However, the method used in
\cite{AzzolliniAandPomponioA-Jmaa,RuizD-JFA} does not apply to
(\ref{eq:1.1}) when $g(x)\not\equiv 0$. In this paper, by
introducing a suitable approximation problem, we try to use the
Theorem 1.1 of \cite{Jeanjean-PRSESM} to get a special (PS) sequence
for $I_\lambda$ based on the weak solutions of the approximation
problem, then to show that this  special (PS) sequence converges to
a solution of problem (\ref{eq:1.1}) in the case of $g(x)\not\equiv
0$.
 We should mention that this kind of idea has been
used in \cite{AmbrosettiAntonioRuizDavid-CCM} to get multiple
solutions to (\ref{eq:1.1}) in the case of $g(x)\equiv 0$ and $p\in (2,5)$.
However, when $g(x)\not\equiv 0$ we cannot prove the boundedness of a (PS) sequence by following
similar idea as those of Lemma 2.6 in \cite{AmbrosettiAntonioRuizDavid-CCM}, here we have to use
an indirect method to do that, see our proof of Lemma \ref{Le3.2}.

However when $p\in (1,2]$, we note that (\ref{eq:1.1})
has no any positive energy solution for $\lambda>0$ large enough (see
Theorem \ref{th:4.1}). Based on this observation, by using the cut-off technique as in
\cite{JeanjeanCoz-AdvDE}(see also
\cite{AzzolliniAveniaPomponio-AIHP,KikuchiHiroaki-ANonSt}) and combining some
delicate analysis, we finally get a positive energy solution for problem (\ref{eq:1.1}) with
$\lambda>0$ small.

\noindent {\bf Notations:} Throughout this paper, we denote
  the standard norms of $H^1(\mathbb{R}^3)$ and $L^p(\mathbb{R}^3)$ by $||\cdot||$ and
$|\cdot|_p$, respectively.

For $p\in [1,5]$, by Sobolev embedding theorem, we have
\begin{equation}\label{embed}
\inf\limits_{|u|_{p+1}=1}\|u\|\triangleq S>0.
\end{equation}
Let
\begin{equation}\label{eq:1.5}
C_p=\frac{(p-1)}{2p}\left(\frac{(p+1)S^{p+1}}{2p}\right)^{1/(p-1)}.
\end{equation}

Our main results are as follows:
\begin{theorem}\label{th:1.1}
If $p\in(2,5)$ and $g(x) \in C^1(\mathbb{R}^3)\cap
L^2(\mathbb{R}^3)$ is a nonnegative function satisfying
\begin{description}
\item $\rm(G1)$: $g(x)=g(|x|)\not\equiv 0$.
\item $\rm(G2)$:  $\langle\nabla g(x),x\rangle\in
L^2(\mathbb{R}^3)$.
\item $\rm(G3)$: $|g|_2<C_p$, where $C_p$ given by
(\ref{eq:1.5}).
\end{description}
Then, for all $\lambda>0$, problem (\ref{eq:1.1}) has at least two nontrivial solutions $\tilde{u}_0$
and $\tilde{u}_1$ such that
$I_\lambda(\tilde{u}_0)<0<I_\lambda(\tilde{u}_1)$.
\end{theorem}

\begin{theorem}\label{th:1.2}
If $p\in(1,2]$, $g(x)$ satisfies $\rm{ (G1)}$  {and}  $\rm{(G3)}$.  Then, only for $\lambda>0$ small,
problem (\ref{eq:1.1}) has two nontrivial solutions $\tilde{u}_0$
and $\tilde{u}_1$ with property
$I_\lambda(\tilde{u}_0)<0<I_\lambda(\tilde{u}_1)$. For $\lambda>0$  large enough, problem (\ref{eq:1.1})  has no any solution with positive
energy.
\end{theorem}

\section {A weak solution with negative energy}
The aim of this section is to get a weak solution with negative energy to problem
(\ref{eq:1.1}), for any $\lambda>0$ and $p\in (1,5)$. With the aid of Ekeland's variational principle, this weak solution is obtained
by seeking a local minimum of the energy functional $I_\lambda$.

\begin{lemma}\label{le2.1}
Let $p\in (1,5)$ and $|g|_2< C_p$ with $C_p$ given by
(\ref{eq:1.5}). Then for the energy functional $I_\lambda$ defined
by (\ref{eq:1.3}), there exist $\alpha>0$ and $\rho >0$ such that
$$I_\lambda(u)\geq\rho>0, \text{\ for \ all \ } \lambda>0 \text{\ and \ } \| u \|=\alpha.$$
\end{lemma}
\textbf{Proof}: For all $\lambda>0$ and $u\in H^1(\mathbb{R}^3)$, by
 Sobolev embedding theorem, we have
\begin{eqnarray}\label{eq:2.1}
I_\lambda(u)&\geq&
\frac{1}{2}\|u\|^2-\frac{1}{(p+1)A_p}\|u\|^{p+1}-|g|_2\|u\|
\nonumber \\
&=&\|u\|(\frac{1}{2}\|u\|-\frac{1}{(p+1)A_p}\|u\|^p-|g|_2),
\end{eqnarray}
where $S_1$ and $A_p$ are given by (\ref{embed}) and (\ref{eq:1.5}).

Set $h(t)=\frac{1}{2}t-\frac{1}{(p+1)A_p}t^p$ for $t\geq 0$. By
direct calculations, we see that
\begin{equation*}
\max\limits_{t\geq 0}h(t)=h(\alpha)= C_p,
\end{equation*}
where $\alpha=[(p+1)S^{p+1}/2p]^{1/(p-1)}$. Then it follows from
(\ref{eq:2.1}) that if $|g|_2< C_p$, there exists
$\rho=\alpha(h(\alpha)-|g|_2)>0$ such that
$I_\lambda(u)|_{\|u\|=\alpha}\geq\rho>0$ for all $\lambda>0$.  $\square$ \\

\begin{theorem}\label{th2.1}
If $p\in (1,5)$, $0\leq g(x)=g(|x|)\in
L^2(\mathbb{R}^3) \setminus \{0\}$ and $|g|_2<C_p$, $C_p$ is given by
(\ref{eq:1.5}). Then for any $\lambda>0$, there exists $u_0\in
H_r^1(\mathbb{R}^3)$ such that
\begin{equation}\label{eq:2.2}
I_\lambda(u_0)=\inf\{I_\lambda(u): u\in H_r^1(\mathbb{R}^3) \ {\rm
and} \   \|u\|\leq \alpha\}<0,
\end{equation}
where $\alpha$ is given by Lemma \ref{le2.1}. Moreover, $u_0$ is a
 solution of problem (\ref{eq:1.1}).
\end{theorem}
\textbf{Proof}: Since $g(x)=g(|x|)\in L^2(\mathbb{R}^3)$, $g(x)\geq
0$ and $g(x)\not\equiv 0$, we can choose a function $v\in
H_r^1(\mathbb{R}^3)$ such that $\int g(x)vdx>0$. Then for $t>0$
small enough, we have
\begin{equation}\nonumber
I_\lambda(tv)=\frac{t^2}{2}\int|\nabla v|^2+ v^2 dx+\frac{\lambda
t^4}{4}\int\phi_v v^2dx-\frac{t^{p+1}}{p+1}\int|v|^{p+1}dx-t\int
g(x)v dx<0.
\end{equation}
This shows that $c_0:=\inf\{I_\lambda(u): u\in \bar{B}_\alpha\}<0$,
where $\bar{B}_\alpha =\{u\in
H_r^1(\mathbb{R}^3):\|u\|\leq\alpha\}$. By Ekeland's variational
principle, there exists $\{u_n\}\subset \bar{B}_\alpha $ such that
$$
{\rm (i)} \  c_0\leq I_\lambda(u_n)\leq c_0+\frac{1}{n},  \ {\rm
and} \  {\rm (ii)} \ I_\lambda(w)\geq
I_\lambda(u_n)-\frac{1}{n}\|w-u_n\| \ {\rm for \ all} \ w\in
\bar{B}_\alpha.
$$
From a standard procedure, see for example \cite{XPZhu-HSZhou}, we
can prove that $\{u_n\}$ is a bounded (PS) sequence of $I_\lambda$.
Then by the compactness of the embedding
$H_r^1(\mathbb{R}^3)\hookrightarrow L^p(\mathbb{R}^3)(2<p<6)$, there
exists $u_0\in H_r^1(\mathbb{R}^3)$ such that
$u_n\stackrel{n}{\rightarrow}u_0$ strongly in $H_r^1(\mathbb{R}^3)$.
Hence $I_\lambda(u_0)=c_0<0$ and $I_\lambda'(u_0)=0$.          $\square$ \\

\section {Positive energy solution for $p \in (2,5)$}

In this section, we aim to prove that problem (\ref{eq:1.1})
has a mountain pass type (positive energy) solution for any $\lambda>0$ and $p\in (2,5)$.
As is known, it is not easy to show that a (PS) sequence of the functional $I_\lambda$
is bounded when $p\in (1,3)$ because of the appearance of nonlocal term of (\ref{eq:1.2}).
 In particular, $p \in (1,2]$ is the hardest case, which we will be deal with in the following section.
To
show the boundedness of a (PS) sequence of $I_\lambda$ in the case of $p \in (2,5)$, it is also nontrivial. Here we have to use a theorem
of \cite{Jeanjean-PRSESM}, which is essentially based on Struwe's monotonicity trick
\cite{Struwe-CMH} and it has been successfully used to handle many
homogeneous elliptic problems, for examples,
\cite{AmbrosettiAntonioRuizDavid-CCM,AzzolliniAveniaPomponio-AIHP}
and the references therein. Motivated by these papers, we apply this
theorem to solve our inhomogeneous elliptic
problem (\ref{eq:1.1}). Let us recall the abstract theorem.
\begin{theorem} \cite[Theorem 1.1]{Jeanjean-PRSESM}\label{th:3.1}
Let $(X,\|\cdot\|)$ be a Banach space, $J\subset\mathbb{R}^+$ an
interval and $(I_\mu)_{\mu\in J}$ a family of $C^1$-functionals on
$X$ of the form
$$I_\mu(u)=A(u)-\mu B(u), \ \forall \ \mu\in J,$$
where $B(u)\geq0$, $\forall \ u\in X$ and $B(u)\rightarrow+\infty$
or $A(u)\rightarrow+\infty$ as $\|u\|\rightarrow\infty$. Assume that
there are two points $v_1,v_2\in X$ such that
$$c(\mu):=\inf_{\gamma\in \Gamma}\max_{t\in[0,1]}I_{\mu}(\gamma(t))>\max\{I_\mu(v_1),I_\mu(v_2)\}, \ {\rm for} \ \mu\in J,$$
where
$$\Gamma=\{\gamma\in
C([0,1],X):\gamma(0)=v_1,\gamma(1)=v_2\}.$$ Then, for almost every
$\mu\in J$, there is a sequence $\{v_n\}\subset X$ such that\\ {\rm\bf (i)}
$\{v_n\}$ is bounded, {\rm\bf(ii)} $I_\mu(v_n)\rightarrow c(\mu)$, {\rm\bf(iii)}
$I'_\mu(v_n)\rightarrow0$ in the dual of $X$.\\
\end{theorem}
In order to applying Theorem \ref{th:3.1} to get a solution to our problem
(\ref{eq:3.1}), we introduce, for any fixed$\lambda>0$, the following approximation
problem
\begin{equation}\label{eq:3.1}
\left\{\begin{array}{ll}
 -\Delta u +  u+\lambda\phi (x) u =\mu|u|^{p-1}u+g(x), \ x\in \mathbb{R}^3, \\
 -\Delta\phi = u^2, \ x\in \mathbb{R}^3,
 \end{array}\right.
\end{equation}
where $\mu\in[1/2,1]$, $p\in(2,5)$ and $g(x)=g(|x|)\in
L^2(\mathbb{R}^3)$.

Let $X=H_r^1(\mathbb{R}^3)$ and $J=[1/2,1]$, and
define $I_{\lambda,\mu}:X\to \mathbb{R}$ by

$$I_{\lambda,\mu}(u)=A(u)-\mu B(u), \rm{with}$$
$$A(u)=\frac{1}{2}\int|\nabla u|^2+ u^2
dx+\frac{\lambda}{4}\int\phi_uu^2dx-\int g(x)udx, \
B(u)=\frac{1}{p+1}\int|u|^{p+1}dx.$$ Then $(I_{\lambda,\mu})_{\mu\in
J}$ is a family of $C^1$-functionals on $X$, $B(u)\geq0$, $\forall \
u\in X$ and $A(u)\geq
\frac{1}{2}\|u\|^2-|g|_2\|u\|\rightarrow+\infty$ as
$\|u\|\rightarrow\infty$.

\begin{lemma}\label{Le3.1}
Let $\lambda>0$ be fixed. Assume that $p\in(2,5)$, $0\leq
(\not\equiv)g(x)=g(|x|)\in L^2(\mathbb{R}^3)$ and $|g|_2<C_p$ with
$C_p$ given by (\ref{eq:1.5}), then
\begin{description}
\item $\rm(i)$ There exist $a, b>0$ and
$e\in H_r^1(\mathbb{R}^3)$ with $\|e\|>b$ such that
\begin{equation}\nonumber
I_{\lambda,\mu}(u)\geq a>0 \ \text{ with }\|u\|=b\ \ \text{ and } \
I_{\lambda,\mu}(e)<0\ \ \text{ for all }\mu\in [1/2,1].
\end{equation}
\item $\rm(ii)$ For any $\mu\in[1/2,1]$,  we have
$$c_{\lambda,\mu}:=\inf_{\gamma\in \Gamma}\max_{t\in[0,1]}I_{\lambda,\mu}(\gamma(t))
>\max\{I_{\lambda,\mu}(0),I_{\lambda,\mu}(e)\},$$
where $\Gamma=\{\gamma\in
C([0,1],H_r^1(\mathbb{R}^3)):\gamma(0)=0,\gamma(1)=e\}$.
\end{description}
\end{lemma}
\textbf{Proof}: (i) Since $I_{\lambda,\mu}(u)\geq I_{\lambda,1}(u)$
for all $u\in H_r^1(\mathbb{R}^3)$ and $\mu\in [1/2,1]$, by Lemma
\ref{le2.1} there exist $a, b>0$, which are independent of $\mu\in
[1/2,1]$, such that $I_{\lambda,1}(u)\geq a>0 \ \text{ with
}\|u\|=b$.

We choose a function $w\in H_r^1(\mathbb{R}^3)\geq(\not\equiv)0$.
Setting $w_t(x)=t^2w(tx)$ for $t>0$, then we have for all $\mu\in
[1/2,1]$,
\begin{eqnarray*}
I_{\lambda,\mu}(w_t)&\leq& \frac{1}{2}\int|\nabla w_t|^2+ w_t^2
dx+\frac{\lambda}{4}\int\phi_{w_t}w_t^2dx-\frac{1}{2(p+1)}\int|w_t|^{p+1}dx
\nonumber \\
&=&\frac{t^3}{2}\int|\nabla w|^2dx+ \frac{t}{2}\int w^2
dx+\frac{\lambda t^3}{4}\int\phi_{w}w^2dx-\frac{
t^{2p-1}}{2(p+1)}\int|w|^{p+1}dx. \nonumber \\
\end{eqnarray*}
Noting that $p\in (2,5)$,  there exists $t_0>0$ large enough, which
is independent of $\mu\in [1/2,1]$,  such that
$I_{\lambda,\mu}(w_{t_0})<0$ for all $\mu\in [1/2,1]$. Hence, (i)
holds by taking $e=w_{t_0}$.

(ii) By the definition of $c_{\lambda,\mu}$, we have for all $\mu\in
[1/2,1]$,
$$c_{\lambda,\mu}\geq c_{\lambda,1}\geq a>0,$$
where $a>0$ is given in (i). Since $I_{\lambda,\mu}(0)=0$ and
$I_{\lambda,\mu}(e)<0$ for all $\mu\in [1/2,1]$, we see that (ii)
holds.   $\square$ \\

By Lemma \ref{Le3.1} and Theorem \ref{th:3.1}, there exists
$\{\mu_j\}\subset [1/2,1]$ such that
\begin{description}
\item (i) $\mu_j\to 1$ as $j\to +\infty$, and
\item (ii) $I_{\lambda,\mu_j}$ has a bounded (PS) sequence $\{u_n^j\}$ at
the level $c_{\lambda,\mu_j}$.
\end{description}
Since the embedding $H_r^1(\mathbb{R}^3)\hookrightarrow
L^p(\mathbb{R}^3)(2<p<6)$ is compact, we can show that for each
$j\in \mathbb{N}$, there exists $u_j\in H_r^1(\mathbb{R}^3)$ such
that $u_n^j\stackrel{n}{\rightarrow}u_j$ strongly in
$H_r^1(\mathbb{R}^3)$ and $u_j$ is a  solution of problem
(\ref{eq:3.1}) with $\mu=\mu_j$. Moreover, we have
\begin{equation}\label{eq:3.2}
0<a\leq I_{\lambda,\mu_j}(u_j)=c_{\lambda,\mu_j}\leq
c_{\lambda,\frac{1}{2}} \ \text{and} \ I_{\lambda,\mu_j}'(u_j)=0, \ \text{for all} \ j\in \mathbb{N}.
\end{equation}
Under the conditions of Theorem \ref{th:1.1}, following the argument
in \cite{BerestyckiLions-ARMA}, we can prove that $u_j$ satisfies
the following type of Pohozaev identity
\begin{equation}\label{eq:3.3}
\int\frac{1}{2}|\nabla u_j|^2+\frac{3}{2}u_j^2
+\frac{5}{4}\lambda\phi_{u_j}u_j^2dx=\int\frac{3\mu_j}{p+1}|u_j|^{p+1}+\left(3g(x)+\langle
x,\nabla g(x)\rangle\right)u_jdx.
\end{equation}
In what follows, we turn to showing that $\{u_j\}$ converges to a solution of problem (\ref{eq:1.1}). For this purpose,  it is necessary to prove that $\{u_j\}$ is  bounded in
$H_r^1(\mathbb{R}^3)$. If $g(x) \equiv 0$, this can be done directly by solving the system of linear equations (\ref{eq:3.2}) and (\ref{eq:3.3}) for $\{|u_j|_2\}$ and $\{|\nabla u_j|_2\}$. However, if $g(x) \not \equiv 0$, this method seems not work well. Here we introduce a new system based on (\ref{eq:3.2}) and (\ref{eq:3.3}), then argue by contradiction.

\begin{lemma}\label{Le3.2}
Under the conditions of Theorem \ref{th:1.1}, if $p\in (2,5)$, then
$\{u_j\}$ is  bounded in $H_r^1(\mathbb{R}^3)$.
\end{lemma}
\textbf{Proof}: We prove the lemma by the following two steps.

{\bf Step 1.} $\{|u_j|_2\}$ is  bounded. \\
By contradiction, we assume that
$|u_j|_2\stackrel{j}{\rightarrow}+\infty$. Let
$v_j=\frac{u_j}{|u_j|_2}$, $X_j=\int |\nabla v_j|^2dx$,
$Y_j=\lambda|u_j|_2^2\int\phi_{ v_j}v_j^2dx$ and
$Z_j=\mu_j|v_j|_{p+1}^{p+1}|u_j|_2^{p-1}$.
  It follows from (\ref{eq:3.2}) that
\begin{equation}\label{eq:3.4}
\left\{\begin{array}{lll} \int\frac{1}{2}|\nabla u_j|^2+\frac{1}{2}
u_j^2+\frac{\lambda}{4}\phi_{u_j}u_j^2dx-\frac{\mu_j}{p+1}\int|u_j|^{p+1}dx-\int
g(x)u_jdx=c_{\lambda,\mu_j},\\
 \int|\nabla u_j|^2+u_j^2
+\lambda\phi_{u_j}u_j^2dx-\int g(x)u_jdx=\mu_j\int|u_j|^{p+1},
 \end{array}\right.
\end{equation}
and
$\{c_{\lambda,\mu_j}\}$ is bounded. Note that  $g(x),
\langle\nabla g(x),x\rangle\in L^2(\mathbb{R}^3)$. Multiplying (\ref{eq:3.3}) and
(\ref{eq:3.4}) by $\frac{1}{|u_j|_2^2}$, we see that
\begin{equation}\label{eq:3.5}
\left\{\begin{array}{lll}
\frac{1}{2}X_j+\frac{1}{4}Y_j-\frac{1}{p+1}Z_j=-\frac{1}{2}+o(1),\\
 \frac{1}{2}X_j
+\frac{5}{4}Y_j-\frac{3}{p+1}Z_j=-\frac{3}{2}+o(1), \\
 X_j+Y_j-Z_j=-1+o(1),
 \end{array}\right.
\end{equation}
where $o(1)$ denotes the quantity tends to zero as
$j\rightarrow+\infty$. For $p\in (2,5)$, solving(\ref{eq:3.5}) we
have
$$X_j=\frac{p-1}{2(2-p)}+o(1).$$
This is a contradiction for $j$ large enough since $X_j\geq 0$ for
all $j\in \mathbb{N}$. Thus, for $p\in (2,5)$, $\{| u_j|_2\}$ is
bounded.

{\bf Step 2.} $\{|\nabla u_j|_2\}$ is  bounded. \\
Similar to the proof of Step 1, we assume by contradiction that
$|\nabla u_j|_2\stackrel{j}{\rightarrow}+\infty$. Let
$w_j=\frac{u_j}{|\nabla u_j|_2}$, $M_j=\lambda|\nabla
u_j|_2^2\int\phi_{w_j}w_j^2dx$ and
$N_j=\mu_j|w_j|_{p+1}^{p+1}|\nabla u_j|_2^{p-1}$,  then multiplying
(\ref{eq:3.4}) by $\frac{1}{|\nabla u_j|_2^2}$ and noting that $\{|
u_j|_2\}$ is bounded, we get
\begin{equation}\label{eq:3.6}
\left\{\begin{array}{lll}
\frac{1}{4}M_j-\frac{1}{p+1}N_j=-\frac{1}{2}+o(1),\\
\frac{5}{4}M_j-\frac{3}{p+1}N_j=-\frac{1}{2}+o(1), \\
M_j-N_j=-1+o(1).
 \end{array}\right.
\end{equation}
From the first and second equations of (\ref{eq:3.6}), we have
$$M_j=2+o(1), \  \  N_j=p+1+o(1).$$
This and the third equation of
(\ref{eq:3.6}) implies that $p=2+o(1)$. So, if $p\neq 2$, we see
that (\ref{eq:3.6}) is impossible. Thus,
for $p\in (2,5)$, $\{|\nabla u_j|_2\}$ is  bounded.  $\square$ \\

\textbf{Proof of Theorem \ref{th:1.1}}: By Lemma \ref{Le3.2} we
can show that $\{u_j\}$ is a bounded (PS) sequence of $I_\lambda$.
Then by the compactness of the embedding
$H_r^1(\mathbb{R}^3)\hookrightarrow L^{p+1}(\mathbb{R}^3)(2<p<5)$, it
follows that for any $\lambda>0$, problem (\ref{eq:1.1}) has a
solution $u_1$ satisfying $I_\lambda(u_1)>0$. Thus, combining Theorem
\ref{th2.1}, we complete the proof of Theorem \ref{th:1.1}.
$\square$ \\

\section {Positive energy solution for $p\in (1, 2]$}
In this section, we claim first that problem (\ref{eq:1.1}) with
$1<p\leq 2$ has no any solution with positive energy for $\lambda>0$ large enough.
\begin{theorem}\label{th:4.1}
Assume that $p\in(1,2]$ and $g(x)\in L^2(\mathbb{R}^3)$(may not be
radially symmetric). Then problem (\ref{eq:1.1}) has no any solution with
positive energy if $\lambda>0$ is large enough.
\end{theorem}
\textbf{Proof:} Let $w\in H^1(\mathbb{R}^3)$ be a solution of
problem (\ref{eq:1.1}). Then $\langle I_\lambda'(w),w\rangle=0$ and
\begin{equation}\label{eq:4.1}
I_\lambda(w)=-\left\{\frac{1}{2}\int|\nabla w|^2+ w^2
dx+\frac{3\lambda}{4}\int\phi_ww^2dx -
\frac{p}{p+1}\int|w|^{p+1}dx\right\}.
\end{equation}
By (20) of \cite{RuizD-JFA}, we have
\begin{equation}\label{eq:4.2}
\sqrt{\lambda/8}\int|w|^3dx\leq\frac{1}{4}\int |\nabla
w|^2dx+\frac{\lambda}{8}\int \phi_{w} w^2dx.
\end{equation}
For $p\in (1, 2]$ and
$\lambda>0$ large enough, it follows from (\ref{eq:4.1}) and  (\ref{eq:4.2})  that
\begin{equation}\nonumber
I_\lambda(w)\leq-\left\{\int\frac{1}{2} w^2+\sqrt{\lambda/2}|w|^3
-\frac{p}{p+1}|w|^{p+1}dx\right\}<0.
\end{equation}
Hence, problem (\ref{eq:1.1}) must have no any solution with positive energy
if $\lambda>0$ is large enough.
$\square$ \\

When $p\in(1,2]$, Theorem \ref{th:4.1} implies that we
may find a solution with
positive energy to  problem (\ref{eq:1.1}) only for $\lambda>0$ small. In this case,  to get a bounded
$\rm (PS)_c  (c>0)$ sequence of $I_\lambda$, following
\cite{JeanjeanCoz-AdvDE}  we introduce the cut-off function $\eta\in
C^\infty(\mathbb{R}^+,\mathbb{R}^+)$ satisfying
\begin{eqnarray}\label{eq:4.3}
\left\{\begin{array}{ll} \eta(t)= 1, \text{ for } t\in[0,1],\\
0\leq\eta(t)\leq 1, \text{ for } t\in (1,2), \\
\eta(t)=0, \text{ for } t\in[2,+\infty), \\
|\eta'|_{\infty}\leq 2,
 \end{array}\right.
\end{eqnarray}
and consider the modified functional
$I_{\lambda,M}:H_r^1(\mathbb{R}^3)\to \mathbb{R}$ defined by
\begin{equation}\label{eq:4.4}
I_{\lambda,M}(u)=\frac{1}{2}\int|\nabla u|^2+ u^2
dx+\frac{\lambda}{4}\int\psi_M(u)\phi_uu^2dx-\frac{1}{p+1}\int|u|^{p+1}dx-\int
g(x)udx,
\end{equation}
where $\psi_M(u)=\eta(\frac{\|u\|^2}{M^2})$ for $M>0$.   If
$g(x)=g(|x|)\in L^2(\mathbb{R}^3)$ and $p\in [1,5]$, we have
$I_{\lambda,M}\in C^1(H_r^1(\mathbb{R}^3),\mathbb{R})$ for each
$\lambda,M>0$ and
\begin{eqnarray}\label{eq:4.5}
\langle I_{\lambda,M}'(u),\varphi\rangle&=&\int\nabla u\nabla
\varphi+ u\varphi
dx+\lambda\psi_M(u)\int\phi_u u\varphi dx   \nonumber \\
&+&
\frac{\lambda}{2}\eta'(\frac{\|u\|^2}{M^2})\frac{1}{M^2}\int\phi_u
u^2dx\int\nabla u\nabla \varphi+ u\varphi dx \nonumber \\
&-&\int|u|^{p-1}u\varphi dx-\int g(x)\varphi dx,
\end{eqnarray}
for any $\varphi\in H_r^1(\mathbb{R}^3)$.

\begin{lemma}\label{Le:4.1}
Assume that $p\in (1,5)$, $0\leq (\not\equiv)g(x)=g(|x|)\in
L^2(\mathbb{R}^3)$ and $|g|_2<C_p$ with $C_p$ given by
(\ref{eq:1.5}). Then the functional $I_{\lambda,M}$ satisfies \\
$\rm(i)$ $I_{\lambda,M}|_{\|u\|=\alpha}>\rho>0$ for all
$\lambda,M>0$,  \\
$\rm(ii)$ For each $M>0$, there exists $e_M\in H_r^1(\mathbb{R}^3)$
with $\|e_M\|>\alpha$ such that $I_{\lambda,M}(e_M)<0$ for all
$\lambda>0$, \\
 where $\alpha,\rho$ are given in Lemma \ref{le2.1}.
 \end{lemma}
\textbf{Proof:} (i) The proof is similar to that of Lemma
\ref{le2.1}. \\

(ii) We choose a function $0\leq v_1\in H_r^1(\mathbb{R}^3)$ such
that $\|v_1\|=1$. By (\ref{eq:4.3}) and (\ref{eq:4.4}),  for each
$M>0$, there exists $t_M\geq 2M>0$ large enough such that
$\psi_M(t_Mv_1)=0$ and $I_{\lambda,M}(t_Mv_1)<0$. Hence, (ii) holds
by taking $e_M=t_Mv_1$.
$\square$ \\

Define
$$c_{\lambda,M}=\inf\limits_{\gamma\in\Gamma_{\lambda,M}}\max\limits_{t\in[0,1]}
I_{\lambda,M}(\gamma(t)),$$ where $\Gamma_{\lambda,M}:=\{\gamma\in
C([0,1],H^1_r(\mathbb{R}^3)):\gamma(0)=0, \gamma(1)=e_M\}$. Then by
Lemma \ref{Le:4.1}, we have
\begin{equation}\label{eq:4.6}
c_{\lambda,M}\geq \rho>0 \ \text{for all} \ \lambda,M>0.
\end{equation}
Furthermore, applying mountain pass theorem, there exists
$\{u_{\lambda,M}^n\}\subset H^1_r(\mathbb{R}^3)$ such that
\begin{equation}\label{eq:4.7}
I_{\lambda,M}(u_{\lambda,M}^n)\stackrel{n}{\rightarrow}c_{\lambda,M}
\ \ \text{and} \ \
(1+\|u_{\lambda,M}^n\|)\|I_{\lambda,M}'(u_{\lambda,M}^n)\|_{H_r^{-1}}\stackrel{n}\rightarrow
0.
\end{equation}
where $H_r^{-1}$ denotes the dual space of $H^1_r(\mathbb{R}^3)$.

\begin{lemma}\label{Le:4.2}
Under the conditions of Lemma \ref{Le:4.1}, let
$\{u_{\lambda,M}^n\}$ be given by (\ref{eq:4.7}), then there exists
$M_0>0$ such that
 $$\limsup\limits_{n\rightarrow+\infty}\|u_{\lambda,M_0}^n\|\leq M_0/2, \ \text{ for all} \ 0<\lambda<M_0^{-3}.$$
\end{lemma}
\textbf{Proof:} Motivated by \cite{KikuchiHiroaki-ANonSt}, we prove the lemma by contradiction. Assume that, for every $M>0$,
there exists $\lambda_M\in (0,M^{-3})$ such that
\begin{equation}\label{eq:4.8}
\limsup\limits_{n\rightarrow+\infty}\|u_{\lambda_M,M}^n\|> M/2.
\end{equation}
For simplicity, we denote $u_{\lambda_M,M}^n$ by $u_n$. By
(\ref{eq:4.8}) and up to a subsequence, we get $\|u_n\|\geq M/2$ for
all $n\in \mathbb{N}$.

From  (\ref{eq:4.4}) and (\ref{eq:4.5}), we have
\begin{eqnarray*}
&&I_{\lambda_M,M}(u_n)-\frac{1}{p+1}\langle
I_{\lambda_M,M}'(u_n),u_n\rangle \nonumber \\
&=&(\frac{1}{2}-\frac{1}{p+1})\|u_n\|^2
+(\frac{1}{4}-\frac{1}{p+1})\lambda_M\int\psi_M(u_n)\phi_{u_n}u_n^2dx
\nonumber \\
&-&\frac{\lambda_M}{2(p+1)}\eta'(\frac{\|u_n\|^2}{M^2})\frac{\|u_n\|^2}{M^2}\int\phi_{u_n}u_n^2dx
-\frac{p}{p+1}\int g(x)u_ndx.
\end{eqnarray*}
Therefore,
\begin{eqnarray}\label{eq:4.9}
&&(\frac{1}{2}-\frac{1}{p+1})\|u_n\|^2-\frac{1}{p+1}\|I_{\lambda_M,M}'(u_n)\|_{H_r^{-1}}\|u_n\| \nonumber \\
&\leq& (\frac{1}{2}-\frac{1}{p+1})\|u_n\|^2 + \frac{1}{p+1}\langle
I_{\lambda_M,M}'(u_n),u_n\rangle \nonumber \\
&=&
I_{\lambda_M,M}(u_n)+\frac{3-p}{4(p+1)}\lambda_M\int\psi_M(u_n)\phi_{u_n}u_n^2dx \nonumber \\
&+&\frac{\lambda_M}{2(p+1)}\eta'(\frac{\|u_n\|^2}{M^2})\frac{\|u_n\|^2}{M^2}\int\phi_{u_n}u_n^2dx
+\frac{p}{p+1}\int g(x)u_ndx.
\end{eqnarray}
By (\ref{eq:4.7}), we have
$I_{\lambda_M,M}(u_n)=c_{\lambda_M,M}+o(1)$, here $o(1)$ denotes the
quantity tends to zero as $n\to +\infty$.

We claim that there exist $M_1,C_1,D_1>0$ such that
\begin{equation}\label{eq:4.10}
c_{\lambda_M,M}\leq C_1\lambda_M M^4+D_1, \ \text{for all} \ M\geq
M_1.
\end{equation}
Let  $v_1$ be the function taken in the proof of Lemma \ref{Le:4.1} (ii), by (\ref{eq:4.4}) we
have
\begin{equation}\label{eq:4.11}
I_{\lambda_M,M}(2Mv_1)\leq 2M^2-\frac{2^{p+1}}{p+1}
|v_1|_{p+1}^{p+1}M^{p+1}.
\end{equation}
Then there exists $M_1>0$ such that $I_{\lambda_M,M}(2Mv_1)<0$ for
all $M\geq M_1$. Thus,
\begin{equation}\label{eq:4.12}
c_{\lambda_M,M}\leq\max\limits_{t\in [0,1]}I_{\lambda_M,M}(t2Mv_1),
\ \text{for all} \ M\geq M_1.
\end{equation}
By (\ref{eq:4.4}) we have
\begin{eqnarray}\label{eq:4.13}
&&\max\limits_{t\in [0,1]}I_{\lambda_M,M}(t2Mv_1)\nonumber \\
&\leq&
\max\limits_{t\in [0,1]}\{2(Mt)^2-\frac{2^{p+1}}{p+1}
|v_1|_{p+1}^{p+1}(Mt)^{p+1}\} %\nonumber \\
 +\max\limits_{t\in
[0,1]}\{\frac{\lambda_M}{4}(2tM)^4\int\phi_{v_1}v_1^2dx\} \nonumber \\
&\leq&\max\limits_{s\geq 0}\{2s^2-\frac{2^{p+1}}{p+1}
|v_1|_{p+1}^{p+1}s^{p+1}\}+C_1\lambda_M M^4 \nonumber \\
&=&D_1+C_1\lambda_M M^4.
\end{eqnarray}
It follows from (\ref{eq:4.12}) and (\ref{eq:4.13}) that
(\ref{eq:4.10}) holds.

By the inequality $\int\phi_uu^2dx\leq C\|u\|^4$ for any $u\in
H^1(\mathbb{R}^3)$ (see \cite{RuizD-JFA}) and noting that
$\psi_M(u_n)=0$ for $\|u_n\|^2\geq 2M^2$, it is easy to see that
\begin{equation}\label{eq:4.14}
\int\psi_M(u_n)\phi_{u_n}u_n^2dx\leq CM^4,
\end{equation}
and
\begin{equation}\label{eq:4.15}
\eta'(\frac{\|u_n\|^2}{M^2})\frac{\|u_n\|^2}{M^2}\int\phi_{u_n}u_n^2dx\leq
CM^4.
\end{equation}

Combining (\ref{eq:4.7}), (\ref{eq:4.9}), (\ref{eq:4.10})
(\ref{eq:4.14}) and (\ref{eq:4.15}), we get for all $M\geq M_1$,
\begin{equation}\label{eq:4.16}
(\frac{1}{2}-\frac{1}{p+1})\|u_n\|^2\leq C_2\lambda_M
M^4+D_2+\frac{p}{p+1}\int g(x)u_ndx,
\end{equation}
where $C_2,D_2>0$ independent of $M$. Then using the inequality
$\int g(x)u_ndx\leq \epsilon\|u_n\|^2+C(\epsilon,|g|_2)$ for any
$\epsilon>0$ and (\ref{eq:4.16}), we deduce that  there exist
$C_3,D_3>0$ independent of $M$ such that for all $ M\geq M_1$,
\begin{equation}\label{eq:4.17}
\|u_n\|^2\leq C_3\lambda_M M^4+D_3.
\end{equation}
Since $\lambda_M\leq M^{-3}$ and $\|u_n\|\geq \frac{M}{2}$,
(\ref{eq:4.17}) is impossible for $M>0$ large enough. Thus we
complete the proof of  this lemma. $\square$ \\

\textbf{Proof of Theorem \ref{th:1.2}}: By Lemma \ref{Le:4.2}
and (\ref{eq:4.3})-(\ref{eq:4.5}), we see that
$\{u_{\lambda,M_0}^n\}\subset H_r^1(\mathbb{R}^3)$ is a bounded (PS)
sequence of $I_\lambda$ for all $0<\lambda<M_0^{-3}$. Moreover, from
(\ref{eq:4.6}) and (\ref{eq:4.7}) we have
$$I_\lambda(u_{\lambda,M_0}^n)=I_{\lambda,M_0}(u_{\lambda,M_0}^n)\stackrel{n}{\rightarrow}c_{\lambda,M_0}\geq
\rho>0.$$ Since
 the embedding $H_r^1(\mathbb{R}^3)\hookrightarrow
L^{p+1}(\mathbb{R}^3)(1<p<5)$ is compact, it follows that for any
$0<\lambda<M_0^{-3}$, problem (\ref{eq:1.1}) has a solution
$\tilde{u}_1$ satisfying $I_\lambda(\tilde{u}_1)>0$. Then by
Theorems \ref{th2.1} and \ref{th:4.1} we complete the proof of
Theorem \ref{th:1.2}.
$\square$ \\

\end{document}